\documentclass[final]{siamart0216}

\usepackage{lipsum}
\usepackage{amsfonts}
\usepackage{amsmath}
\usepackage{graphicx}
\usepackage{epstopdf}
\usepackage{algorithmic}
\usepackage{listings}
\usepackage{pgf,tikz}
\usepackage{multirow}
\usepackage{upgreek}
\usepackage{bm}

\usetikzlibrary{arrows}
\usetikzlibrary{patterns}
\ifpdf
  \DeclareGraphicsExtensions{.eps,.pdf,.png,.jpg}
\else
  \DeclareGraphicsExtensions{.eps}
\fi

\newcommand{\bn}{\ensuremath{{\bm{\nu}}}}
\newcommand{\bx}{\ensuremath{{\mathbf{x}}}}
\newcommand{\by}{\ensuremath{{\mathbf{y}}}}

\newcommand{\bp}{\ensuremath{{\mathbf{p}}}}
\newcommand{\avg}[1]{\{\!\{{#1}\}\!\}}
\newcommand{\e}{\mathrm{e}}
\newcommand{\ii}{\mathrm{i}}
\newcommand{\dx}[1][x]{\,\mathrm{d}#1}

\newcommand{\interior}{^\text{--}}
\newcommand{\exterior}{^\text{+}}

\newcommand{\mat}[1]{#1}           % mat = matrix (from discretised operator)
\newcommand{\bop}[1]{\mathsf{#1}}  % bop = boundary operator
\newcommand{\pop}[1]{\mathcal{#1}} % pop = potential operator

\newcommand{\TheTitle}{Product algebras for Galerkin discretisations of boundary integral operators and their applications} 
\newcommand{\TheAuthors}{T. Betcke, M. W. Scroggs, and W. \'{S}migaj}

\headers{Product Algebras for Galerkin discretisations}{\TheAuthors}

\title{{\TheTitle}\thanks{This work was funded by EPSRC Grants EP/I030042/1 and EP/K03829X/1}}

\author{
  Timo Betcke\thanks{Department of Mathematics, University College London, UK
    (\email{t.betcke@ucl.ac.uk}).}
  \and
  Matthew W. Scroggs\thanks{Department of Mathematics, University College London, UK
  (\email{matthew.scroggs.14@ucl.ac.uk}).}
  \and
  Wojciech \'{S}migaj\thanks{Simpleware Ltd., Exeter, UK (\email{w.smigaj@simpleware.com}).}
}

\usepackage{amsopn}

\ifpdf
\hypersetup{
  pdftitle={\TheTitle},
  pdfauthor={\TheAuthors}
}
\fi

\begin{document}

\maketitle

% REQUIRED
\begin{abstract}
  Operator products occur naturally in a range of regularized boundary integral equation formulations. However, while a Galerkin discretisation only depends on the domain space and the test (or dual) space of the operator, products require a notion of the range. In the boundary element software package Bempp we have implemented a complete operator algebra that depends on knowledge of the domain, range and test space. The aim was to develop a way of working with Galerkin operators in boundary element software that is as close to working with the strong form on paper as possible while hiding the complexities of Galerkin discretisations. In this paper, we demonstrate the implementation of this operator algebra and show, using various Laplace and Helmholtz example problems, how it significantly simplifies the definition and solution of a wide range of typical boundary integral equation problems.
\end{abstract}

% REQUIRED

% REQUIRED

\section{Introduction}
\label{sec:introduction}
A typical abstract operator problem can be formulated as
$$
\bop{A}\phi = f,
$$
where $\bop{A}$ is an operator mapping from a Hilbert space $\mathcal{H}_1$ into another Hilbert space $\mathcal{H}_2$ with the unknown $\phi\in\mathcal{H}_1$ and known $f\in\mathcal{H}_2$. Many modern operator preconditioning strategies depend on the idea of having a regulariser $\bop{R}:\mathcal{H}_2\to\mathcal{H}_1$ and solving the equation
\begin{equation}
\label{eq:regularized_eq}
\bop{R}\bop{A}\phi = \bop{R}f
\end{equation}
instead. This is particularly common in the area of boundary integral equations, where integral operators can be efficiently preconditioned by operators of opposite order. Now suppose that we want to discretise \eqref{eq:regularized_eq} using a standard Galerkin method. The discretised problem is
\begin{equation}
\label{eq:discretized_eq}
\mat{R}M^{-1}\mat{A}\bm{\phi} = \mat{R}\bm{f},
\end{equation}
where $\mat{R}$ and $\mat{A}$ are the Galerkin discretisations of $\bop{R}$ and $\bop{A}$, respectively, and $\bm{f}$ is the vector of coefficients of the projection of $f$ onto the finite dimensional subspace of $\mathcal{H}_2$. The matrix $\mat{M}$ is the mass matrix between the basis functions of the finite dimensional subspaces of $\mathcal{H}_2$ and $\mathcal{H}_1$.

In order to solve \eqref{eq:discretized_eq}, we have to assemble all involved matrices, form the right-hand side, implement a function that evaluates $\mat{R}\mat{M}^{-1}\mat{A}\bm{v}$ for a given vector $\bm{v}$, and then solve \eqref{eq:discretized_eq} with GMRES or another iterative solver of choice.
Ideally, we would not have to deal with these implementational details and just directly write the following code.
\begin{lstlisting}[language=Python]
A = operator(...)
R = operator(...)
f = function(...)
phi = gmres(R * A, R * f)
\end{lstlisting}
Note that at the end, the solution \verb@phi@ is again a function
object. In order for this code snippet to work and the mass matrix $M$ to be assembled automatically, either the
implementation of the operator product needs to be aware of the test
space of $\bop{A}$ and domain space of $\bop{R}$, or the software definition of
these operators need to contain information about their ranges.
In this paper we will follow this second approach by defining the notion of the strong form of a Galerkin discretisation and demonstrate its benefits.

An implementation of a product algebra based on this idea is contained in the Python/C++ based boundary element library Bempp (\url{www.bempp.com}) \cite{Smigaj2015}, originally developed by the authors of this paper. Bempp is a comprehensive library for the solution of boundary integral equations for Laplace, Helmholtz and Maxwell problems. The leading design principle of Bempp is to allow a description of BEM problems in Python code that is as close to the mathematical formulation as possible, while hiding implementational details of the underlying Galerkin discretisations. This allows us to formulate complex block operator systems such as those arising in Calder\'{o}n preconditioned formulations of transmission problems in just a few lines of code. Initial steps towards a Bempp operator algebra were briefly described in \cite{Smigaj2015} as part of a general library overview. The examples in this paper are based on the current version (Bempp 3.3), which has undergone significant development since then and now contains a complete and mature product algebra for operators and grid functions.

As examples for the use of an operator algebra in more complex settings, we discuss: the efficient assembly of the hypersingular operator via a representation using single layer operators; the assembly of Calder\'on projectors and the computation of their spectral properties and the Calder\'{o}n preconditioned solution of acoustic transmission problems.

A particular challenge is the design of product algebras for Maxwell problems. The stable discretisation of the electric and magnetic field operators for Maxwell problems requires the use of a non-standard skew symmetric bilinear form. The Maxwell case is discussed in much more detail in \cite{Scroggs2017}.

The paper is organised as follows. In Section \ref{sec:integral_operators} we review basic definitions of boundary integral operators for Laplace and Helmholtz problems. In Section \ref{sec:product_algebra} we introduce the basic concepts of a Galerkin product algebra and discuss some implementational details. Section \ref{sec:fast_assembly} then gives a first application to the fast assembly of hypersingular operators for Laplace and Helmholtz problems. Then, in Section \ref{sec:block_systems} we discuss block operator systems at the example of Calder\'on preconditioned transmission problems. The paper concludes with a summary in Section \ref{sec:conclusions}.

While most of the mathematics presented in this paper is well known among specialists, the focus of this paper is on hiding mathematical complexity of Galerkin discretisations. With the wider penetration and acceptance of high-level scripting languages such us Matlab, Python and Julia in the scientific computing community, we now have the tools and structures to make complex computational operations accessible for a wide audience of non-specialist users, making possible the fast dissemination of new algorithms and techniques beyond traditional mathematical communities.

\section{Boundary integral operators for scalar Laplace and Helmholtz problems, and their Galerkin discretisation}
\label{sec:integral_operators}

In this section, we give the basic definitions of boundary integral operators for Laplace and Helmholtz problems and some of their properties needed later. More detailed information can be found in e.g. \cite{Steinbach08, Sau11}. 

We assume that $\Omega\subset\mathbb{R}^3$ is a piecewise smooth bounded Lipschitz domain with boundary $\Gamma$. By $\Omega\exterior := \mathbb{R}^3\backslash\overline{\Omega}$ we denote the exterior of $\Omega$. We denote by $\gamma_0^{\pm}$ the associated interior (-)  and exterior (+) trace operators and by $\gamma_1^{\pm}$ the interior and exterior normal derivative operators. We always assume that the normal direction $\bn$ points outwards into $\Omega\exterior $.

The average of the interior and exterior trace is defined as $\avg{\gamma_0 f} := \frac{1}{2}\left(\gamma_0\exterior f + \gamma_0\interior f\right)$. Correspondingly, the average normal derivative is defined as $\avg{\gamma_1 f} := \frac{1}{2}\left(\gamma_1\exterior f + \gamma_1\interior f\right)$.

\subsection{Operator definitions}
\label{sec:op_definitions}

We consider a function $\phi\interior\in H^1(\Omega)$ satisfying the Helmholtz equation $-\Delta \phi\interior - k^2\phi\interior =0$, where $k\in\mathbb{R}$. By Green's representation theorem we have
\begin{align}
\label{eq:int_green}
\phi\interior(\bx) &= \left[\pop{V}\gamma_1\interior{\phi\interior}\right](\bx)-\left[\pop{K}\gamma_0\interior\phi\interior\right](\bx),&\bx\in\Omega
\end{align}
for the single layer potential operator $\mathcal{V}:{H^{-1/2}(\Gamma)}\to H_\text{loc}^{1}(\Omega\cup\Omega\exterior )$ defined by
\begin{align*}
\left[\pop{V}\mu\right](\mathbf{\bx})&=\int_{\Gamma}G(\bx,\by)\mu(\by)\dx[s(\by)],&\mu\in H^{-1/2}(\Gamma)
\end{align*}
and the double layer potential operator $\mathcal{K}:{H^{1/2}(\Gamma)}\to H_\text{loc}^{1}(\Omega\cup\Omega\exterior )$ defined by
\begin{align*}
\left[\pop{K}\xi\right](\mathbf{\bx})&=\int_{\Gamma}\frac{\partial G(\bx,\by)}{\partial \bn(\by)}\xi(\by)\dx[s(\by)],&\xi\in H^{1/2}(\Gamma).
\end{align*}
Here, $G(\bx, \by) := \frac{\e^{\ii k|\bx-\by|}}{4\uppi|\bx-\by|}$ is the associated Green's function. If $k=0$, we obtain the special case of the Laplace equation $-\Delta u=0$.

We now define the following boundary operators as the average of the interior and exterior traces of the single layer and double layer potential operators:
\begin{itemize}
\item The \textit{single layer} boundary operator $\bop{V}:H^{-1/2}(\Gamma)\to H^{1/2}(\Gamma)$ defined by
\begin{align*}
\left[\bop{V}\mu\right](\bx) &= \avg{\gamma_0\pop{V}\mu}(\bx),&\mu\in H^{-1/2}(\Gamma),~\bx\in\Gamma.
\end{align*}
\item The \textit{double layer} boundary operator $\bop{K}:H^{1/2}(\Gamma)\to H^{1/2}(\Gamma)$ defined by
\begin{align*}
\left[\bop{K}\xi\right](\bx) &= \avg{\gamma_0\pop{K}\xi}(\bx),&\xi\in H^{1/2}(\Gamma),~\bx\in\Gamma.
\end{align*}
\item The \textit{adjoint double layer} boundary operator $\bop{K}':H^{-1/2}(\Gamma)\to H^{-1/2}(\Gamma)$ defined by
\begin{align*}
\left[\bop{K}'\mu\right](\bx) &= \avg{\gamma_1\pop{V}\mu}(\bx),&\mu\in H^{-1/2}(\Gamma),~\bx\in\Gamma.
\end{align*}
\item The \textit{hypersingular} boundary operator $\bop{W}:H^{1/2}(\Gamma)\to H^{-1/2}(\Gamma)$ defined by
\begin{align*}
\left[\bop{W}\xi\right](\bx) &= -\avg{\gamma_1\pop{K}\xi}(\bx),&\xi\in H^{1/2}(\Gamma),~\bx\in\Gamma.
\end{align*}
\end{itemize}
Applying the interior traces $\gamma_0\interior$ and $\gamma_1\interior$ to the Green's representation formula \eqref{eq:int_green}, and taking into account the jump relations of the double layer and adjoint double layer boundary operators on the boundary $\Gamma$ \cite[Section 6.3 and 6.4]{Steinbach08}
we arrive at
\begin{equation}
\label{eq:calderon_int}
\begin{bmatrix}\gamma_0\interior\phi\interior\\ \gamma_1\interior\phi\interior\end{bmatrix}
= \left(\tfrac{1}{2}\bop{Id}+\bop{A}\right)\begin{bmatrix}\gamma_0\interior\phi\interior\\ \gamma_1\interior\phi\interior\end{bmatrix}
\end{equation}
with
\begin{equation}
\label{eq:multitrace_operator}
\bop{A} := \begin{bmatrix}
-\bop{K} & \bop{V}\\
\bop{W} & \bop{K}'
\end{bmatrix},
\end{equation}
which holds almost everywhere on $\Gamma$. The operator $\mathcal{C}\interior := \frac{1}{2}\bop{Id} + \bop{A}$ is also called the interior Calder\'{o}n projector. If $\phi\exterior $ is a solution of the exterior Helmholtz equation $-\Delta\phi\exterior -k^2\phi\exterior  = 0$ in $\Omega\exterior $ with boundary condition at infinity
$$
\lim_{|\bx|\to\infty} |\bx|\left(\frac{\partial}{\partial|\bx|}\phi\exterior -\ii k\phi\exterior \right)=0
$$
for $k\neq 0$
and
$$
\lim_{|\bx|\to\infty}|\phi\exterior (\bx)| = \mathcal{O}\left(\frac{1}{|\bx|}\right)
$$
for $k=0$, Green's representation formula is given as
\begin{equation}
\label{eq:ext_green}
\phi\exterior (\bx) = \left[\pop{K}\gamma_0\exterior \phi\exterior \right](\bx)-\left[\pop{V}\gamma_1\exterior {\phi\exterior }\right](\bx),~\bx\in\Omega\exterior .
\end{equation}
Taking the exterior traces $\gamma_0\exterior $ and $\gamma_1\exterior $ now gives the system of equations
\begin{equation}
\label{eq:calderon_ext}
\begin{bmatrix}\gamma_0\exterior \phi\exterior \\ \gamma_1\exterior \phi\exterior \end{bmatrix}
= \left(\tfrac{1}{2}\bop{Id}-\bop{A}\right)\begin{bmatrix}\gamma_0\exterior \phi\exterior \\ \gamma_1\exterior \phi\exterior .\end{bmatrix}
\end{equation}
with associated exterior Calder\'{o}n projector $\mathcal{C}\exterior :=\frac{1}{2}\bop{Id}-\bop{A}$.

\subsection{Galerkin discretisation of integral operators}
\label{sec:galerkin}

Let $\mathcal{T}_h$ be a triangulation of $\Gamma$ with $N$ piecewise flat triangular elements $\tau_j$ and $M$ associated vertices $\bp_i$. We define the function space $S_h^{0}$ of elementwise constant functions $\phi_j$ such that
$$
\phi_j(\bx) =
\begin{cases}
1,&x\in\tau_j\\
0,&\text{otherwise},
\end{cases}
$$
and the space $S_h^{1}$ of globally continuous, piecewise linear hat functions $\rho_i$ such that
$$
\rho_i(\bp_\ell) =
\begin{cases}
1,&i = \ell\\
0,&\text{otherwise}.
\end{cases}
$$
Denote by $\langle u,v\rangle_{\Gamma}$ the standard surface dual form $\int_{\Gamma}u(\bx)\overline{v(\bx)}\dx[s(\bx)]$ of two functions $u$ and $v$. By restricting $H^{1/2}(\Gamma)$ onto $S_h^{1}$ and $H^{-1/2}(\Gamma)$ onto $S_h^0$, we obtain the Galerkin discretizations $\mat{V}$, $\mat{K}$, $\mat{K'}$, $\mat{W}$ defined as
$$
\begin{array}{cc}
\left[\mat{V}\right]_{ij} := \langle \bop{V}\phi_j, \phi_i\rangle_{\Gamma},&
\left[\mat{K}\right]_{ij} := \langle \bop{K}\rho_j, \phi_i\rangle_{\Gamma}\\
\left[\mat{K'}\right]_{ij} := \langle \bop{K}'\phi_j, \rho_i\rangle_{\Gamma},&
\left[\mat{W}\right]_{ij} := \langle \bop{W}\rho_j, \rho_i\rangle_{\Gamma}
\end{array}
$$
From this definition it follows that $\mat{K'} = \mat{K}^T$. A computable expression of $\mat{W}$ using weakly singular integrals is given in Section \ref{sec:fast_assembly}.

A problem with this definition of discretisation spaces is that $S_h^0$ and $S_h^1$ have a different number of basis functions, leading to non-square matrices $K$ and $K'$. Hence, it is only suitable for discretisations of integral equations of the first-kind involving only $\bop{V}$ or $\bop{W}$ on the left-hand side. There are two solutions to this.
\begin{enumerate}
\item Discretise both spaces $H^{1/2}(\Gamma)$ and $H^{-1/2}(\Gamma)$ with the continuous space $S_h^1$. This works well if $\Gamma$ is sufficiently smooth. However, if $\Gamma$ has corners then Neumann data in $H^{-1/2}(\Gamma)$ is not well represented by continuous functions.
\item Instead of the space $S_h^0$ use the space of piecewise constant functions $\phi_{D}$ on the dual grid which is obtained by associating each element of the dual grid with one vertex of the original grid. We denote this piecewise constant space by $S_{D,h}^0$. With this definition of piecewise constant functions also the matrix $\mat{K}$ is square. Moreover, the mass matrix between the basis functions in $S_h^0$ and $S_{D,h}^0$ is inf-sup stable \cite{Hiptmair2006,Buffa2007}. 
\end{enumerate}

\section{Galerkin product algebras and their implementation}
\label{sec:product_algebra}

In this section we discuss the product of Galerkin discretisations of abstract Hilbert space operators and how a corresponding product algebra can be implemented in software. While the mathematical basis is well known, most software libraries do not support a product algebra, making implementations of operator based preconditioners and many other operations more cumbersome than necessary. This section proposes a framework to elegantly support operator product algebras in general application settings. The formalism introduced here is based on Riesz mappings between dual spaces. A nice introduction in the context of Galerkin discretizations is given in \cite{Kirby2010}.

\subsection{Abstract formulation}
\label{sec:abstract_algebra}

Let $\bop{A}:\mathcal{H}_\bop{A}^\text{dom}\to\mathcal{H}_\bop{A}^\text{ran}$ and $\bop{B}:\mathcal{H}_\bop{B}^\text{dom}\to\mathcal{H}_\bop{B}^\text{ran}$ be operators mapping between Hilbert spaces. If $\mathcal{H}_\bop{A}^\text{ran}\subset\mathcal{H}_\bop{B}^\text{dom}$ the product
\begin{equation}
\label{eq:op_product}
g = \bop{B}\bop{A}f
\end{equation}
is well defined in $\mathcal{H}_\bop{B}^\text{ran}$. We now want to evaluate this product using Galerkin discretisations of the operators $\bop{A}$ and $\bop{B}$. 

Let $\mathcal{H}_\bop{A}^\text{dual}$ be dual to $\mathcal{H}_\bop{A}^\text{ran}$ with respect to a given dual pairing $\langle\cdot, \cdot\rangle_\bop{A}:\mathcal{H}_\bop{A}^\text{ran}\times\mathcal{H}_\bop{A}^\text{dual}\to \mathbb{C}$. Correspondingly, we define the space $\mathcal{H}_\bop{B}^\text{dual}$ as dual space to $\mathcal{H}_\bop{B}^\text{ran}$ with respect to a dual pairing $\langle\cdot,\cdot\rangle_\bop{B}$.

Defining the function $q = \bop{A}f$, the operator product \eqref{eq:op_product} can equivalently be written as
\begin{equation}
\begin{array}{r@{}l}
q &{}= \bop{A}f\nonumber\\
g &{}= \bop{B}q\nonumber.
\end{array}
\end{equation}
Rewriting this system in its variational form leads to the problem of finding
$(q, g)\in\mathcal{H}_\bop{A}^\text{ran}\times\mathcal{H}_\bop{B}^\text{ran}$ such that
\begin{equation}
\label{eq:var_product_system}
\begin{array}{r@{}l}
\langle q, \mu\rangle_\bop{A} &{}= \langle\bop{A}f,\mu\rangle_\bop{A}\\
\langle g, \tau\rangle_\bop{B} &{}= \langle\bop{B}q, \tau\rangle_\bop{B}
\end{array}
\end{equation}
for all $(\mu, \tau)\in\mathcal{H}_\bop{A}^\text{dual}\times\mathcal{H}_\bop{B}^\text{dual}$.
We now introduce the finite dimensional subspaces $\mathcal{V}_{h,\bop{X}}^\text{dom}\subset\mathcal{H}_\bop{X}^\text{dom}$, $\mathcal{V}_{h,\bop{X}}^\text{ran}\subset\mathcal{H}_\bop{X}^\text{ran}$ and $\mathcal{V}_{h,\bop{X}}^\text{dom}\subset\mathcal{H}_\bop{X}^\text{dual}$ with basis functions $\zeta_{\bop{X},j}^\text{dom}$, $\zeta_{\bop{X},i}^\text{ran}$, $\zeta_{\bop{X},\ell}^\text{dual}$ for $\bop{X}=\bop{A},\bop{B}$. In what follows we assume that the dimension of $\mathcal{V}_{h,\bop{X}}^\text{dual}$ is identical to the dimension of $\mathcal{V}_{h,\bop{X}}^\text{ran}$ and that the associated dual-pairing is inf-sup stable in the sense that 
$$
\sup_{\xi_{X}^{\text{dual}}\in V_{h,X}^{\text{dual}}}\frac{\langle \xi_X^{\text{ran}}, \xi_X^{\text{dual}}\rangle_{X}}{\|\xi_X^{\text{dual}}\|_{\mathcal{H}_{\bop{X}}^{\text{dual}}}}\geq c_X\|\xi_X^{\text{ran}}\|_{\mathcal{H}_{\bop{X}}^{\text{ran}}},\quad \forall  \xi_{X}^{\text{ran}}\in V_{h,X}^{\text{ran}}
 $$
for some $c_X>0$, implying that the associated mass matrix is invertible. The discrete version of \eqref{eq:var_product_system} is now given as
\begin{eqnarray}
\mat{M}_\bop{A}\bm{q} = \mat{A}\bm{f},\nonumber\\
\mat{M}_\bop{B}\bm{g} = \mat{B}\bm{q}\nonumber,
\end{eqnarray}
where $\left[\mat{M}_\bop{X}\right]_{\ell,i} = \langle \phi_{\bop{X},i}^\text{ran}, \phi_{\bop{X},\ell}^\text{dual}\rangle_\bop{X}$, $\bop{X}=\bop{A},\bop{B}$. are the mass matrices of the dual pairings. The vectors $\bm{f}$, $\bm{q}$ and $\bm{g}$ are the vectors of coefficients of the corresponding functions. Combining both equations we obtain
$$
\bm{q} = \mat{M}_\bop{B}^{-1}\mat{B}\mat{M}_\bop{A}^{-1}\mat{A}\bm{f}.
$$
The matrix $\mat{A}$ is also called the discrete \textit{weak form} of the operator $\bop{A}$. This motivates the following definition.
\begin{definition}Given the discrete weak form $\mat{A}$ defined as above. We define the associated \textit{discrete strong form} as the matrix
$$
\mat{A}^{S} :=\mat{M}_\bop{A}^{-1}\mat{A}.
$$
\end{definition}
Note that $M_A^{-1}$ is the discrete Riesz map from the dual space into the range space of $\mat{A}$ \cite{Kirby2010}.
The notation of the discrete strong form allows us to define a Galerkin product algebra as follows.
\begin{definition}
Given the operator product $\bop{C}:=\bop{B}\bop{A}$. We define the associated discrete operator product weak form as
$$
\mat{C} := \mat{B}\odot\mat{A}:= \mat{B}\cdot\mat{A}^{S} = \mat{B}\mat{M}_\bop{A}^{-1}\mat{A}
$$
and the associated discrete strong form as
$$
\mat{C}^{S} := \mat{M}_\bop{B}^{-1}\left(\mat{B}\odot\mat{A}\right).
$$
\end{definition}
We note that the a direct discretisation $\langle \bop{B}\bop{A}\phi_{\bop{A},j}^\text{dom},\phi_{\bop{B},\ell}^\text{dual}\rangle$ is usually not identical to $\mat{C}$ as the latter is computed as the solution of the operator system \eqref{eq:var_product_system} whose discretisation error also depends on the space $\mathcal{V}_{h,\bop{A}}^\text{ran}$ and the corresponding discrete dual. However, the discretisation of the operator product $\bop{B}\bop{A}$
can rarely be computed directly and solving \eqref{eq:var_product_system} is usually the only possibility to evaluate this product.

This discrete operator algebra is associative since
$$
\left(\mat{C}\odot\mat{B}\right)\odot\mat{A}=
\mat{C}\mat{M}_\bop{B}^{-1}\mat{B}\mat{M}_\bop{A}^{-1}\mat{A}=
\mat{C}\odot\left(\mat{B}\odot\mat{A}\right).
$$
Moreover, since a Hilbert space is self-dual in its natural inner product $(\cdot, \cdot)$ the discretization
$$
\left[\mat{M}_\bop{A}^{dom}\right]_{ij} = \left(\zeta_{\mat{A}, i}^{\text{dom}}, \zeta_{\mat{A}, j}^{\text{dom}}\right)
$$
of the identity operator $\bop{Id}_\bop{A}^\text{dom}$ is the right unit element with respect to this discrete operator algebra. Correspondingly, the matrix $\mat{M}_{\bop{A}}^{\text{ran}}$ is the left unit element.

We have so far considered the approximation of the weak form $\langle \bop{B}\bop{A}\phi_{\bop{A},j}^\text{dom},\phi_{\bop{B},\ell}^\text{dual}\rangle$, where the operator $\bop{B}$ acts on $\bop{A}\phi_{\bop{A},j}^\text{dom}$. However, there are situations where we want a discrete approximation of the
product $\langle \bop{A}\phi_{\bop{A},j}^\text{dom}, \bop{B}\phi_{\bop{B},\ell}^\text{dom}\rangle$ for $\bop{B}:\mathcal{H}_\bop{B}^\text{dom}\to\mathcal{H}_\bop{A}^\text{dual}$. An example for the assembly of hypersingular operators will be given later. Note that if $\by$ is a coefficient vector of a function $\phi\in \mathcal{H}_\bop{B}^\text{dom}$ then $\tilde{\by} = M_B^{-1}\mat{B}\by$ is the coefficient vector to the Galerkin approximation of $\tilde{\phi} = \bop{B}\phi$. Hence, a discrete approximation of the weak form $\langle \bop{A}\phi_{\bop{A},j}^\text{dom}, \bop{B}\phi_{\bop{B},\ell}^\text{dom}\rangle$ is given by
$$
%\mat{C'} = 
\mat{B}^H\cdot\mat{M}_\bop{B}^{-H}\cdot\mat{A} = \left[\mat{B}^{S}\right]^HA.
$$

This motivates the following definition.
\begin{definition} We define the dual discrete product weak form associated with the operators $\bop{A}$ and $\bop{B}$ as
\begin{equation}
\label{eq:left_product}
\mat{B}\odot_D\mat{A} := \mat{B}^H\cdot\mat{M}_\bop{B}^{-H}\cdot\mat{A}
\end{equation}
and the associated discrete strong form as
$$
\mat{C}:=M_\bop{B,A}^{-1}\left(\mat{B}\odot_D\mat{A}\right).
$$
where $M_\bop{B,A}$ is the mass matrix between the domain space of $B$ and the range space of $A$.
\end{definition}

\subsection{Example:~Operator preconditioned Dirichlet problems}
As a first example, we describe the formulation of an operator preconditioned interior Dirichlet problem using the above operator algebra. We want to solve
\begin{align*}
-\Delta \phi\interior -k^2\phi\interior &= 0 &\text{in }&\Omega\\
\gamma_0\phi\interior&=g&\text{on }&\Gamma
\end{align*}
for a given function $g\in H^{1/2}(\Gamma)$. From the first line of \eqref{eq:calderon_int} we obtain that
$$
\gamma_0\interior\phi\interior = \left(\tfrac{1}{2}\bop{Id}-\bop{K}\right)\gamma_0\interior\phi\interior+\bop{V}\gamma_1\interior\phi\interior.
$$
Substituting the boundary condition, we obtain the integral equation of the first kind
\begin{equation}
\label{eq:dirichlet_problem}
\bop{V}\gamma_1\interior\phi\interior = \left(\tfrac{1}{2}\bop{Id}+\bop{K}\right)g.
\end{equation}
The operator $\bop{V}:H^{-1/2}(\Gamma)\to H^{1/2}(\Gamma)$ is a pseudodifferential operator of order $-1$ and can be preconditioned by the hypersingular operator $\bop{W}:H^{1/2}(\Gamma)\to H^{-1/2}(\Gamma)$, which is a pseudodifferential operator of order $1$ \cite{Steinbach98, Hiptmair2006}. We arrive at the preconditioned problem
\begin{equation}
\label{eq:dirichlet_problem_prec}
\bop{W}\bop{V}\gamma_1\interior\phi\interior = \bop{W}\left(\tfrac{1}{2}\bop{Id}+\bop{K}\right)g.
\end{equation}
Note that the operator $\bop{W}$ is singular if $k=0$. In that case a rank-one modification of the hypersingular operator can be used \cite{Steinbach98}.
For the Galerkin discretisation of \eqref{eq:dirichlet_problem_prec} we use the standard $L^2$ based dual pairing $\langle \cdot,\cdot\rangle$ defined by
$$
\langle u, v\rangle = \int_{\Gamma} u(x)\overline{v}(x)dx,\quad u, v\in L^2(\Gamma)
$$
and note that the spaces $H^{1/2}(\Gamma)$ and $H^{-1/2}(\Gamma)$ are dual with respect to this dual pairing.

For the discretisation of the operators we use the spaces $S_{D,h}^{0}$ and $S_h^{1}$ as described in Section \ref{sec:galerkin}. Using the notation introduced in Section \ref{sec:abstract_algebra} we obtain the discrete system
\begin{equation}
\label{eq:discrete_laplace_prec_system}
\mat{W}\odot\mat{V}\bm{x} = \mat{W}\odot\left(\tfrac{1}{2}\mat{M}+\mat{K}\right)\bm{g},
\end{equation}
where $\bm{x}$ is the vector of coefficients of the unknown function $\phi_h$ in the basis $S_{D,h}^{0}$. The matrix $\mat{M}$ is the discretisation of the identity operator on $H^{1/2}(\Gamma)$. If in addition we want to use Riesz (or mass matrix) preconditioning we can simply take the discrete strong forms of the product operators on the left and right hand side of \eqref{eq:discrete_laplace_prec_system}.

In terms of mathematics the definition of the discrete strong form is simply a notational convenience. We could equally write \eqref{eq:discrete_laplace_prec_system} by directly inserting the mass matrix inverses. The main advantage of an operator product algebra first becomes visible in a software implementation that directly supports the notions of discrete strong forms and operator products. This is described below.

\subsection{Basic software implementation of an operator algebra}
Based on the definition of a discrete product algebra for Galerkin discretisations, we can now discuss the software implementation. Two concepts are crucial: namely that of a \textit{grid function}, which represents functions defined on a grid; and that of an \textit{operator}, which maps grid functions from a discrete domain space into a discrete range space.

\subsubsection{Grid functions}

We start with the description of a grid function. A basic grid function object is defined by a discrete function space and a vector of coefficients on the space. However, for practical purposes this is not always sufficient. Consider the following situation of multiplying the discrete single layer operator $\mat{V}$, discretised with the space of piecewise constant functions $S_h^0$ and a vector of coefficients $\bm{f}$.
The result $\bm{y}=\mat{V}\bm{f}$ is defined as $y_i = \sum_{j=1}^nf_j\langle  V\phi_j, \phi_i\rangle$. Since the single layer operator maps onto $H^{1/2}(\Gamma)$ we would like to obtain a suitable vector of coefficients $\tilde{y}$ of piecewise linear functions in $S_h^{1}$ such that
$$
\bm{y} = \mat{M}\tilde{\bm{y}},
$$
where $\mat{M}$ is the rectangular mass-matrix between the spaces $S_h^0$ and $S_h^1$. Solving for $\tilde{\bm{y}}$ is only possible in a least-squares sense. Moreover, for these two spaces the matrix $M$ may even be ill-conditioined or singular in the least-squares sense, making it difficult to obtain a good approximation in the range space. Hence, we also allow the definition of a grid function purely through the vector of coefficients into the dual space.

The constructors to define a grid function either through coefficients in a given space or through projections into a dual space are defined as follows.
\begin{lstlisting}[language=Python]
fun = GridFunction(space, coefficients=...)
fun = GridFunction(space, dual_space=..., projections=...)
\end{lstlisting}
Associated with these two constructors are two methods that extract the vectors of coefficients or projections.
\begin{lstlisting}[language=Python]
coeffs = fun.coefficients()
proj = fun.projections(dual_space)
\end{lstlisting}
If the grid function is initialised with a coefficient vector, then the first operation just returns this vector. The second operation sets up the corresponding mass matrix \verb@M@ and returns the vector \verb@M * coeffs @. If the grid function is initialised with a vector of projections and a corresponding dual space then access to the coefficients results in a solution of a linear system if the space and dual space have the same number of degrees of freedom. Otherwise, an exception is thrown. If the \verb@projections@ method is called and the given dual space is identical to the original dual space on initialisation the vector \verb@projections@ is returned. Otherwise, first a conversion to coefficient form via a call to \verb@coefficients()@ is attempted.

This dual representation of a grid function via either a vector of coefficients or a vector of projections makes it possible to represent functions in many standard situations, where a conversion between coefficients and projections is mathematically not possible and not necessary for the formulation of a problem.

\subsubsection{Operators}

Typically, in finite element discretisation libraries the definition of an operator requires an underlying weak form, a domain space and a test space. However, to support the operator algebra introduced in Section \ref{sec:abstract_algebra} the range space is also required. Hence, we represent a constructor for a boundary operator in the following form.
\begin{lstlisting}[language=Python]
op = operator(domain, range_, dual_to_range, ...)
\end{lstlisting}
Here, the objects \verb@domain@, \verb@range_@ and \verb@dual_to_range@ describe the finite dimensional domain, range and dual spaces. Each operator provides the following two methods.
\begin{lstlisting}[language=Python]
discrete_weak_form = op.weak_form()
discrete_strong_form = op.strong_form()
\end{lstlisting}
The first one returns the standard discrete weak form while the second one returns the discrete strong form. The \verb@discrete_weak_form@ and
\verb@discrete_strong_form@ are objects that implement at least a matrix-vector routine to multiply a vector with the corresponding discrete operator. The multiplication with the inverse of the mass matrix in the strong form is implemented via computing an LU decomposition and solving the associated linear system.

Important for the performance is caching. The weak form is computed in the first call to the \verb@weak_form()@ method and then cached. Correspondingly, the LU decomposition necessary for the strong form is computed only once and then cached.

\subsubsection{Operations on operators and grid functions}

With this framework the multiplication \verb@res_fun = op * fun@ of a boundary operator \verb@op@ with a grid function \verb@fun@ can be elegantly described in the following way:
\begin{lstlisting}[language=Python]
result_fun = GridFunction(
    space=op.range_, 
    dual_space=op.dual_to_range,
    projections=op.weak_form() * fun.coefficients)
\end{lstlisting}
Alternatively, we could have more simply presented the result as
\begin{lstlisting}[language=Python]
result_fun = GridFunction(
    space=op.range_,
    coefficients=op.strong_form() * fun.coefficients)
\end{lstlisting}
However, the latter ignores that there may be no mass matrix transformation
available that could map from the discrete dual space to the discrete range
space.

As an example, we present a small code snippet from Bempp that maps the constant function $f(\bx)=1$ on the boundary of the cube to the function $g=\bop{V}f$, where $\bop{V}$ is the Laplace single layer boundary operator. $f$ is represented in a space of piecewise constant functions on the dual grid and $g$ is represented in a space of continuous, piecewise linear functions, reflecting the smoothing properties of the Laplace single layer boundary operator. The following lines define the cube grid with an element size of $h=0.1$ and the spaces of piecewise constant functions on the dual grid, and continuous, piecewise linear continuous functions on the primal grid. 
\begin{lstlisting}[language=Python]
grid = bempp.api.shapes.cube(h=0.1)
const_space = bempp.api.function_space(grid, "DUAL", 0)
lin_space = bempp.api.function_space(grid, "B-P", 1)
\end{lstlisting}
We would like to remark on the parameter \verb@B-P@ (\textit{barycentric-polynomial}) in the code given above for the function space definitions. Since the piecewise constant functions are defined on the dual grid, we are working with the barycentric refinement of the original grid \cite{Buffa2007}. Hence, the piecewise linear functions on the primal grid also need to be defined over the barycentric refinement (denoted by the parameter \verb@B-P@) as the discretisation routines require the same refinement level for the domain and dual to range space. Mathematically, the standard space of continuous, piecewise linear functions over the primal grid and the space \verb@B-P@ over the barycentric refinement are identical.

We now define the operator and the constant grid function.
For the grid function the coefficient vector is created
via the NumPy routine \verb@ones@, taking as input the number
of degrees of freedom in the space.
\begin{lstlisting}[language=Python]
op = bempp.api.operators.boundary.laplace.single_layer(
    const_space, lin_space, const_space)
fun = bempp.api.GridFunction(
    const_space, 
    coefficients=np.ones(const_space.global_dof_count))
\end{lstlisting}
We can now multiply the operator with the function and plot the result.
\begin{lstlisting}[language=Python]
result = op * fun
result.plot()
\end{lstlisting}
\begin{figure}
\centering
\begin{tabular}{cc}
\includegraphics[width=5cm]{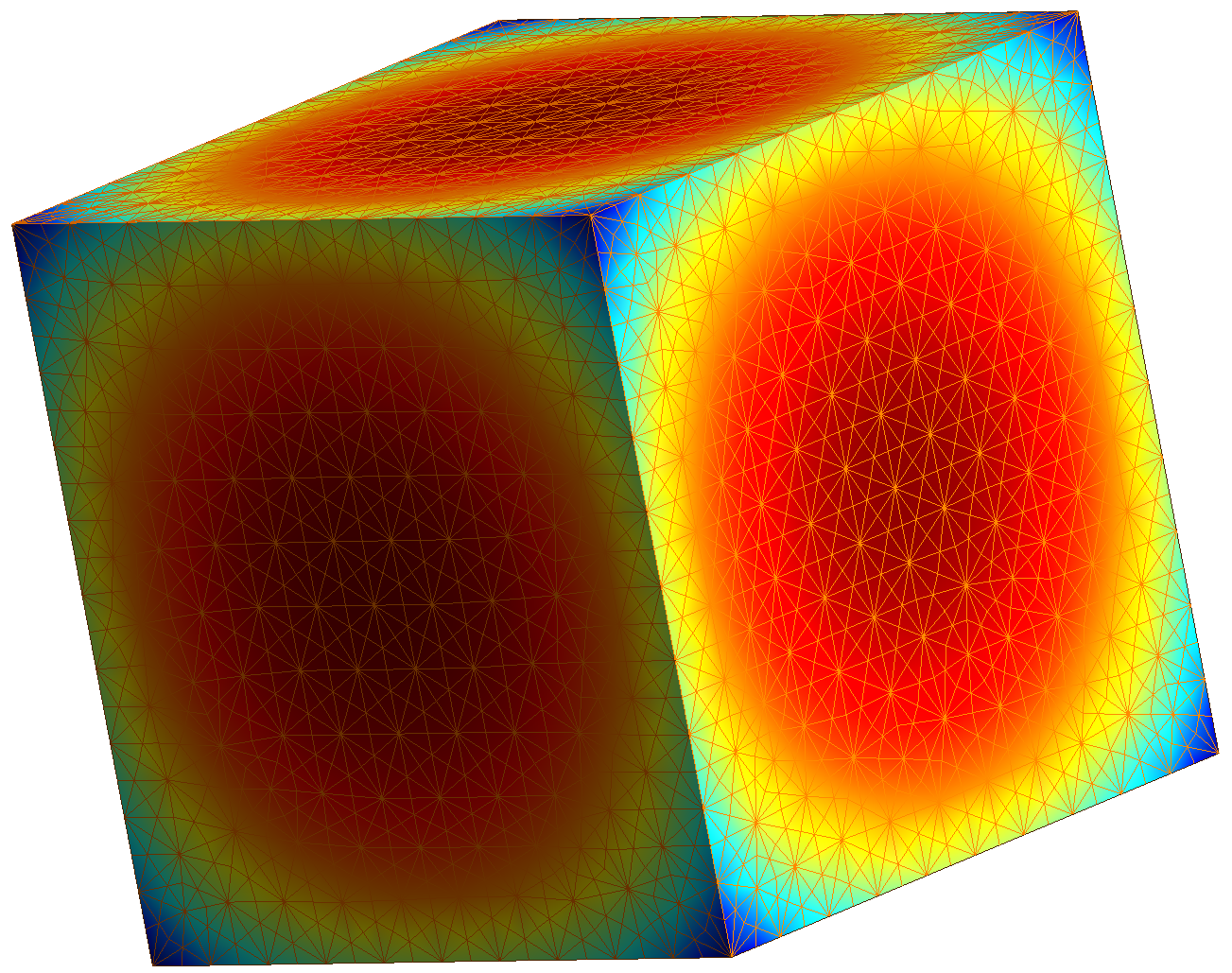} & 
\includegraphics[width=5cm]{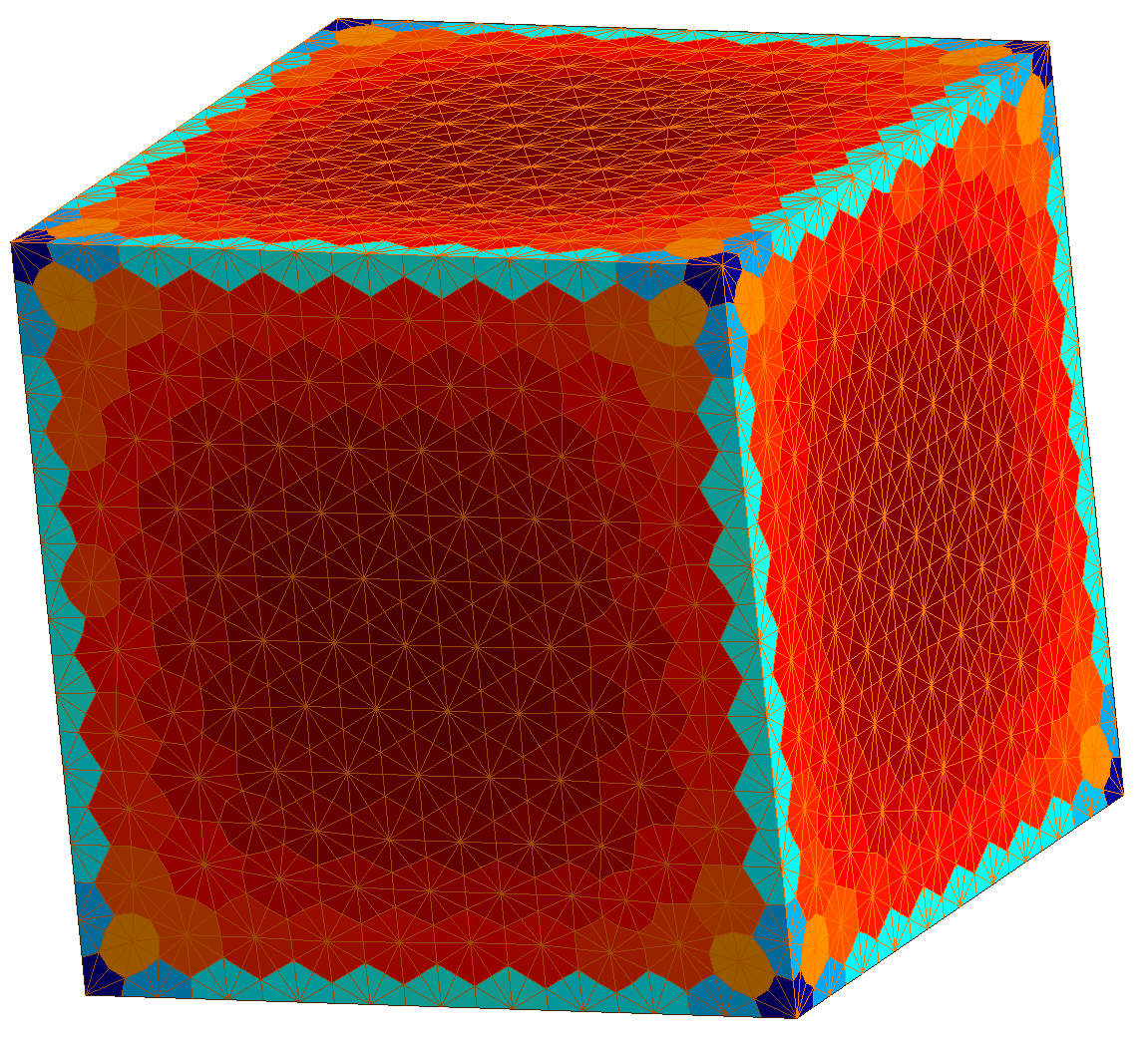}
\end{tabular}
\label{fig:cube}
\caption{Left: The Laplace single layer operator applied to a constant function on the boundary of a cube. Right: The Laplace hypersingular operator applied to the function on the left.}
\end{figure}
The output is the left cube shown in Figure \ref{fig:cube}. It is a continuous function in ${H}^{1/2}(\Gamma)$. The right cube in Figure \ref{fig:cube} shows the result of multiplying the Laplace hypersingular operator defined by
\begin{lstlisting}[language=Python]
op = bempp.api.operators.boundary.laplace.hypersingular(
    lin_space, const_space, lin_space)
\end{lstlisting}
with the function on the left. Since the hypersingular operator maps into $H^{-1/2}(\Gamma)$, the appropriate range space consists of piecewise constant functions, and the result of the discrete operation correspondingly uses a space of piecewise constant functions.

Under the condition that the operations mathematically make sense and operators and functions are correctly defined this mechanism always maps grid function objects into the right spaces under the action of a boundary operator while hiding all the technicalities of Galerkin discretisations.

The internal implementation of the product of two operators is equally simple in this framework. Given two operators \verb@op1@ and \verb@op2@. Internally,
the \verb@weak_form()@ method of the product \verb@op1 * op2@ is defined as follows.
\begin{lstlisting}[language=Python]
def weak_form():
    return op1.weak_form() * op2.strong_form()
\end{lstlisting}
Correspondingly, the strong form of the product is implemented as:
\begin{lstlisting}[language=Python]
def strong_form():
    return op1.strong_form() * op2.strong_form()
\end{lstlisting}
Internally, the product of two discrete operators provides a matrix-vector routine that successively applies the two operators to a given vector. If \verb@op1@ and \verb@op2@ implement caching then an actual discretisation of a weak form is only performed once, and the product of the two operators is performed with almost no overhead.

It is very easy to wrap standard iterative solvers to support this operator algebra. Suppose we want to solve the product system \eqref{eq:discrete_laplace_prec_system}.
Using an operator algebra wrapper to any standard GMRES (such as the one in SciPy \cite{scipy}) the solution to the system \eqref{eq:discrete_laplace_prec_system} now takes the form
\begin{lstlisting}[language=Python]
solution, info = gmres(W * V, W * (.5 * ident + K) * g)
\end{lstlisting}
with \verb@solution@ being a grid function that lives in the correct space of piecewise constant functions.
The definition of such a GMRES routine is as follows:
\begin{lstlisting}[language=Python]
def gmres(A, b, ...):
    from scipy.sparse import linalg
    x, info = linalg.gmres(
        A.weak_form(), 
        b.projections(A.dual_to_range),
        ...)
    return GridFunction(A.domain, coefficients=x), info
\end{lstlisting}    
The product algebra automatically converts \verb@W * V@ into a new object
that provides the correct space attributes and a \verb@weak_form@ method as defined above. Similarly, the right-hand side \verb@b@ is evaluated into a vector
with the \verb@projections@ method. The full Bempp implementation provides among other options also a keyword attribute \verb@use_strong_form@. If this is set to true then inside the GMRES routine the solution is computed as
\begin{lstlisting}[language=Python]
x, info = linalg.gmres(A.strong_form(), b.coefficients)
\end{lstlisting}
This corresponds to standard Riesz (or mass matrix) preconditioning and comes naturally as part of this algebra. Note that we have left out of the description checks that the spaces of the left and right hand side are compatible. In practice, this should be done by the code as sanity check.

%Finally, in order to implement the dual product $\mat{B}\odot_D\mat{A}$ we require the transpose or adjoint (complex conjugate transpose) of an operator. By default, the definition of an operator does not include the dual space of the domain space. But this is required to form the operator adjoint as the range space of the adjoint is just the dual space of the domain space. Hence, the adjoint is defined as
%\begin{lstlisting}[language=Python]
%adjoint_op = op.adjoint(domain_dual)
%\end{lstlisting}
%with \verb@adjoint_op@ being an operator with domain space \verb@op.dual_to_range@, range space \verb@domain_dual@ and dual space \verb@op.domain@. The operator \verb@adjoint_op@ has a \verb@weak_form@ method which returns a discrete operator representing the action of the adjoint of the weak form of \verb@op@ on a given vector.
Finally, the weak form of the dual product $\mat{B}\odot_D\mat{A}$ can be be implemented as
\begin{lstlisting}[language=Python]
def weak_form():
    return B.strong_form().adjoint() * A.strong_form()
\end{lstlisting}
The range space and domain space of the dual product are the same as that of $A$ while the dual space is the same as the domain space of $B$.

\subsection{A note on the performance of the operator algebra}
The operator algebra described above relies on being able to perform fast
mass matrix LU decompositions and solves. In finite element methods LU decompositions with a mass matrix can be as expensive as solves with a stiffness matrix.
In BEM the situation is quite different. Even with the utilisation of fast methods such as FMM (fast multipole method \cite{Greengard1987}) or hierarchical matrices \cite{Hackbusch2015}, the assembly and matrix-vector product of a boundary operator is typically much more expensive than assembling a mass matrix and performing an LU decomposition of it. Therefore, mass matrix operations can be essentially treated as on-the-fly operations compared to the rest. One potential problem is the complexity of the LU decomposition of a mass matrix over a surface function space on $\Gamma$. For banded systems the complexity of Gaussian elimination scales like $\mathcal{O}(n)$. However, a closed surface has a higher element connectivity than a standard plane in 2d and we cannot expect a simple $\mathcal{O}(n)$ scaling even with reordering. In practice though, this has made so far little difference and we have used the SuperLU code provided by SciPy for the LU decomposition and surface linear system solves on medium size BEM problems with hundreds of thousands of surface elements without any noticeable performance issues, and we expect little performance overhead even for very large problems with millions of unknowns as the FMM or hierarchical matrix operations on the operators have much larger effective costs and significantly more complex data structures to operate on.

\section{The fast assembly of hypersingular boundary operators}
\label{sec:fast_assembly}

The weak form of the hypersingular boundary operator can, after integration by parts, be represented as \cite{Hamdi1981, Nedelec1982}
\begin{equation}
\label{eq:hypersing_weak_form}
\begin{array}{r@{}l}
\mat{W}_{ij} &{}=\displaystyle \frac{1}{4\uppi}\int_{\Gamma}\int_{\Gamma}\frac{\e^{\ii k|\bx-\by|}}{|\bx-\by|}\langle\textbf{curl}_{\Gamma}\rho_i(\bx),\textbf{curl}_{\Gamma}\rho_j(\by)\rangle_2\dx[s(\by)]\dx[s(\bx)]\\[3mm]
\displaystyle &{}\displaystyle-\frac{k^2}{4\uppi}\int_{\Gamma}\int_{\Gamma}\frac{\e^{\ii k|\bx-\by|}}{|\bx-\by|}\rho_i(\bx)\rho_j(\by)\langle \bn(\bx),\bn(\by)\rangle_2\dx[s(\by)]\dx[s(\bx)],
\end{array}
\end{equation}
where the basis and test function $\rho_j$ and $\rho_i$ are basis functions in $S_h^1$. Both terms in \eqref{eq:hypersing_weak_form} are now weakly singular and can be numerically evaluated.

% However, \eqref{eq:hypersing_weak_form} is not only important for evaluating the singular near-field of integral operators, but also for the far field when the supports of $\rho_i$ and $\rho_j$ have a nonzero distance. In that case \eqref{eq:hypersing_weak_form} is identical to the direct formula
% $$
% \mat{W}_{ij} = -\int_{\Gamma}\int_{\Gamma}\frac{\partial^2 G(\bx, \by)}{\partial \bn(\bx)\partial\bn(\by)}\rho_i(\bx)\rho_j(\by)ds(\bx)ds(\by).
% $$
% The kernel of this integral is not asymptotically smooth on boundaries with corners due to the two normal derivatives and even on smooth domains is less smooth than the single layer kernel due to the higher order derivatives. This increases the complexity of a low-rank representation of the kernel, making fast methods based on low-rank kernel representations (e.g. FMM or $\mathcal{H}$-matrices) significantly more expensive.

However, \eqref{eq:hypersing_weak_form} motivates another way of assembling the hypersingular operator, which turns out to be significantly more efficient in many cases. In both terms of \eqref{eq:hypersing_weak_form}, a single layer kernel is appearing. We can use this and represent $\mat{W}$
in the form
\begin{equation}
\label{eq:hyp_sum}
\mat{W} = \sum_{j=1}^3\mat{P}_j^T\cdot\mat{V}\cdot\mat{P}_j  -k^2\sum_{j=1}^3\mat{Q}_j^T\cdot\mat{V}\cdot\mat{Q}_j,
\end{equation}
where we now only need to assemble a single layer boundary operator $\mat{V}$ with smooth kernel in a space of discontinuous elementwise linear functions, and the $\mat{P}_j$ and $\mat{Q}_j$ are sparse matrices. $\mat{P}_j$ maps a continuous piecewise linear function to the $j$th component of its surface curl and $\mat{Q}_j$ scales the basis functions with the contributions of $\bn$ in the $j$th component in each element. If $k=0$ (Laplace case) the second term in \eqref{eq:hyp_sum} becomes zero and it would even be sufficient to use a space of piecewise constant functions to represent $\mat{V}$.

This evaluation trick is well known and is suitable for discretising the hypersingular operator with continuous, piecewise linear basis functions on flat triangles. The disadvantage is that an explicit representation of the sparse matrices $P_j$ and $Q_j$ is necessary. This representation depends on the polynomial order and dof numbering of the space implementation.

In the following we use the product algebra concepts to write the representation \eqref{eq:hypersing_weak_form} in a form that generalises to function spaces of arbitrary order on curved triangular elements without requiring details of the dof ordering in the implementation. Given a finite dimensional trial space $V_h^\text{trial}$ with basis $\theta_1,\dots, \theta_L$ and a corresponding test space $V_h^\text{test}$ with basis $\xi_1,\dots,\xi_{L'}$ we define the discrete sparse surface operators
\begin{align}
\left[\mat{C^{\ell}}\right]_{ij} &= \langle \left[\textbf{curl}_{\Gamma}\theta_j\right]_{\ell}, \xi_i\rangle_{\Gamma},\nonumber\\
\left[\mat{N^{\ell}}\right]_{ij} &= \langle \theta_j\left[\bn\right]_{\ell}, \xi_i\rangle_{\Gamma}\nonumber.
\end{align}
The operator $\mat{C^{\ell}}$ weakly maps a function $f$ to its elementwise $\ell$th surface curl component, and the operator $\mat{N^{\ell}}$ weakly multiplies a function $f$ with the $\ell$th component of the surface normal direction.

We can now represent the hypersingular operator as
\begin{equation}
\label{eq:hyp_compound_representation}
\mat{W} = \sum_{j=1}^3 \mat{C^{j}}\odot_D\mat{V}\odot\mat{C^{j}} -k^2\sum_{j=1}^3\mat{N^{j}}\odot_D\mat{V}\odot\mat{N^{j}}. 
\end{equation}
The dual multiplication $\odot_D$ in \eqref{eq:hyp_compound_representation} acts on the test functions and the right multiplication $\odot$ acts on the trial functions. Let $V_h^{m,\text{cont}}$ be a globally continuous, elementwise polynomial function space of order $m$ and denote by $V_h^{m,\text{disc}}$ the corresponding space of discontinuous elementwise polynomial functions of order $m$. Then the operators in \eqref{eq:hypersing_weak_form} have the following domain, range and dual spaces.
\vspace{\baselineskip}
\begin{center}
\begin{tabular}{c|c|c|c}
Operator & domain & range & dual \\
\hline
$\mat{W}$ & $V_h^{m,\text{cont}}$ & $V_h^{m,\text{disc}}$ & $V_h^{m,\text{cont}}$\\
$\mat{V}$ & $V_h^{m,\text{disc}}$ & $V_h^{m,\text{disc}}$ & $V_h^{m,\text{disc}}$\\
$\mat{N^{j}}$ & $V_h^{m,\text{cont}}$ & $V_h^{m,\text{disc}}$ & $V_h^{m,\text{disc}}$\\
$\mat{C^{j}}$ & $V_h^{m,\text{cont}}$ & $V_h^{m,\text{disc}}$ & $V_h^{m,\text{disc}}$
\end{tabular}
\end{center}
\vspace{\baselineskip}
We note that \eqref{eq:hypersing_weak_form} only requires inverses of dual parings on $V_h^{m,\text{disc}}$ with itself as dual space and not dual pairings between $V_h^{m,\text{disc}}$ and $V_h^{m,\text{cont}}$ which are not invertible. If $k=0$ we can use spaces of order $m-1$ for $V$ and the dual and range space of $C$ since then the second sum in \eqref{eq:hyp_sum} vanishes and the first sum only contains products of derivatives of the basis and trial functions. Also, we have chosen the discontinuous function space $V_h^{m,\text{disc}}$ as range space of $\mat{V}$. This guarantees that the result in \eqref{eq:hyp_compound_representation} has the correct range space.

In terms of standard matrix products \eqref{eq:hyp_compound_representation} has the form
$$
\mat{W} = \sum_{j=1}^3 \left[\mat{C^j}\right]^T\cdot \mat{M}^{-T} \cdot \mat{V}\cdot\mat{M}^{-1}\cdot\mat{C^j} -k^2\sum_{j=1}^3\left[\mat{N^j}\right]^T\cdot \mat{M}^{-T} \cdot \mat{V}\cdot\mat{M}^{-1}\cdot\mat{N^j},
$$
where $M$ is the mass matrix associated with the space $V_h^{m, disc}$ of discontinuous basis functions. Hence, $M$ is elementwise block-diagonal and therefore $M^{-1}$ is too, and we can efficiently directly compute $M^{-1}$ as a sparse matrix. We can then accumulate the sparse matrix products in the sum above to obtain \eqref{eq:hyp_sum} with $P_j = \mat{M}^{-1}\cdot \mat{C^j}$ and $Q_j = \mat{M}^{-1}\cdot \mat{N^j}$. In Bempp the whole implementation of the hypersingular operator can be written as follows.
\begin{lstlisting}[language=Python]
D = ZeroBoundaryOperator(...)
for i in range(3):
    D += C[i].dual_product(V) * C[i]
    D += -k**2 * N[i].dual_product(V) * N[i]
\end{lstlisting}
Due to efficient caching strategies, all operators, including the mass matrices and their inverses, are computed only once. Hence, there is minimal overhead from using a high-level expressive formulation.

In Figure \ref{fig:hyp_assembly_modes}, we compare times and memory requirements for the hierarchical matrix assembly of the hypersingular boundary operator on the unit sphere with wavenumber $k=1$ using basis functions in $S_h^1$. The left column shows the standard assembly based on \eqref{eq:hypersing_weak_form} and $S_h^1$ basis functions. The middle column shows results for assembling the operator directly on a larger space of piecewise linear discontinuous functions using the weak form \eqref{eq:hypersing_weak_form} and then projecting down to basis functions in $S_h^1$, that is $\mat{W}=\mat{P}^T\mat{W}_\text{disc}\mat{P}$ for a sparse matrix $P$ that maps from $S_h^1$ to a space of piecewise linear discontinuous functions. This assembly allows matrix compression directly on the elementwise basis functions instead of only compressing on nodal basis functions after summing up the elementwise contributions. However, in the case of the hypersingular operator, this leads to larger memory consumption through the larger matrix size on the discontinuous space, but not faster assembly times. The interesting case is the single layer formulation in \eqref{eq:hyp_compound_representation}. Even though the single layer operator is assembled on the larger discontinuous space it compresses better since it is a smoothing operator and therefore leads to around twice as fast assembly times. The price is a larger memory size compared to the standard assembly. If this is not of concern then the single layer based assembly is preferrable. Note that the evaluation of the matrix-vector product using \eqref{eq:hyp_compound_representation} requires six multiplications with the single layer operator. So if a large number of matrix-vector products is needed this can become a bottleneck.

\begin{figure}
\center
\begin{tabular}{c|c|c|c|c|c|c}
\multirow{2}{*}{N (cont/discont)} & \multicolumn{2}{|c|}{Standard} & \multicolumn{2}{|c|}{projection} & \multicolumn{2}{|c}{via single layer}\\
& time & mem & time & mem & time & mem \\
\hline
258 / 1536     & 0.7 s &1.0 MiB    &  1.0 s &  31 MiB   &  0.4 s    & 15.5 MiB          \\
1026 / 6144     & 3.9 s &9.0 MiB    &  3.5 s    &  177 MiB   &  1.6 s    &   88 MiB       \\
4098 / 24576    & 19.6 s &60.1 MiB    &  15.2 s    &  907 MiB   &  7.7 s    &   467 MiB       \\
16386 / 98304    & 1.6 m & 345 MiB    &  1.4 m    &  4.4 GiB   &  39.4 s   &   2.2 GiB      \\
65538 / 393216    & 7.6 m & 1.79 GiB    &  8.6 m    &  21.5 GiB   &  3.9 m    &    11.0 GiB 
\end{tabular}
\caption{Time and memory for the assembly of the hypersingular operator using the standard weak form on the continuous space, discontinuous assembly with projection spaces or a single layer formulation. In the latter two cases only the assembly time and memory of the boundary operator is given. Assembly time and memory requirements for the sparse operators are negligble. }
\label{fig:hyp_assembly_modes}
\end{figure}

\section{Block operator systems}
\label{sec:block_systems}

Block operator systems occur naturally in boundary element computations since we are typically dealing with pairs of corresponding Dirichlet and Neumann data whose relationship is given by the Calder\'{on} projector shown in \eqref{eq:calderon_int} for the interior problem and \eqref{eq:calderon_ext} for the exterior problem. In this section we want to demonstrate some interesting computations with the Calder\'{o}n projector which can be very intuitively performed in the framework of block operator extensions of the product algebra.

Within the Bempp framework, a blocked operator of given block dimension $(m, n)$ is defined as

\begin{lstlisting}[language=Python]
blocked_operator = bempp.api.BlockedOperator(m, n)
\end{lstlisting}

We can now assign individual operators to the blocked operator by e.g.
\begin{lstlisting}[language=Python]
blocked_operator[0, 1] = laplace.single_layer(...)
\end{lstlisting}

Not every entry of a blocked operator needs to be assigned a boundary operator. Empty positions are
automatically treated as zero operators. However, we require the following conditions before
computations with blocked operators can be performed:
\begin{itemize}
\item There can be no empty rows or columns of the blocked operator.
\item All operators in a given row must have the same \verb@range@ and \verb@dual_to_range@ space.
\item All operators in a given column must have the same \verb@domain@ space.
\end{itemize}

These conditions are easily checked while assigning components to a blocked operator. The weak form of a blocked operator is obtained as
\begin{lstlisting}[language=Python]
discrete_blocked_operator = blocked_operator.weak_form()
\end{lstlisting}

This returns an operator which performs a matrix-vector product by splitting up the input vector into its components with respect to the columns of the blocked operator, performs multiplications with the weak forms of the individual components, and then assembles the result vector back together again.

The interesting case is the definition of a strong form. Naively, we could just take the strong forms of the individual component operators. However, since each strong form involves the solution of a linear system with a mass matrix we want to avoid this. Instead, we multiply the discrete weak form of the operator from the left with a block diagonal matrix whose block diagonal components contain the inverse mass matrices that map from the dual space in the corresponding row to the range space. This works due to the compatibility condition that all test and range spaces within a row must be identical.

\subsection{Stable discretisations of Calder\'{o}n projectors}

With the concept of a block operator we now have a simple framework to work with Calder\'{o}n projectors $\mathcal{C}^{\pm} = \left(\frac12\bop{Id} \mp \bop{A}\right)$ with $\bop{A}$ defined as in \eqref{eq:multitrace_operator}. For the sake of simplicity in the following we use the Calder\'{o}n projector $\mathcal{C}\exterior $ for the exterior problem. The interior Calder\'{o}n projector $\mathcal{C}\interior$ is treated in the same way. Remember that both operators are defined on the product space $H^{1/2}(\Gamma) \times H^{-1/2}(\Gamma)$

Two properties are fundamental to Calder\'{o}n projectors. First, $\left(\mathcal{C}\exterior \right)^2 = \mathcal{C}\exterior $; and second, if $U = \begin{bmatrix}\gamma_0\exterior u,~\gamma_1\exterior u\end{bmatrix}^T$ is the Cauchy data of an exterior Helmholtz solution $u$ satisfying the Sommerfeld radiation condition, it holds that $U=\mathcal{C}\exterior U$, or equivalently $\mathcal{C}\interior U = 0$.

Based on the product algebra framework introduced in this paper we can easily represent these properties on a discrete level to obtain a numerical Calder\'{o}n projector up to the discretisation error.

As an example, we consider the Calder\'{o}n projector on the unit cube with wavenumber $k=2$. Assembling the projector within the Bempp product operator framework is simple, and corresponding functions are already provided.
\begin{lstlisting}[language=Python]
k = 2
from bempp.api.operators.boundary.sparse \
    import multitrace_identity
from bempp.api.operators.boundary.helmholtz \
    import multitrace_operator
calderon = .5 * multitrace_identity(grid, spaces='dual') \
           - multitrace_operator(grid, k, spaces='dual')
\end{lstlisting}
In this code snippet, the option \verb@spaces='dual'@ automatically discretises
the Calder\'{o}n projector using stable dual pairings of continuous, piecewise linear
spaces on the primal grid, and piecewise constant functions on the dual grid.

\begin{figure}
\center
\includegraphics[width=10cm]{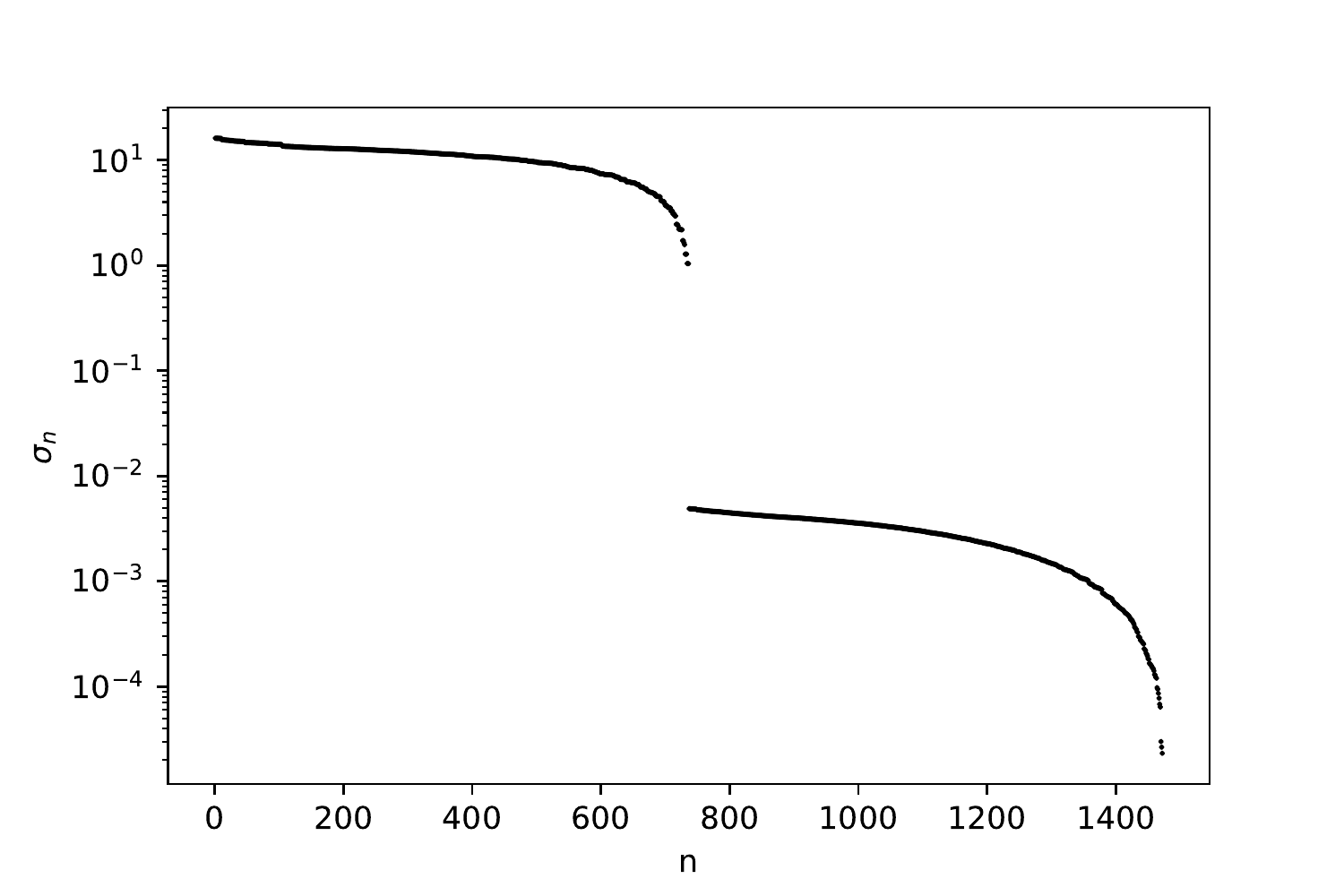}
\includegraphics[width=10cm]{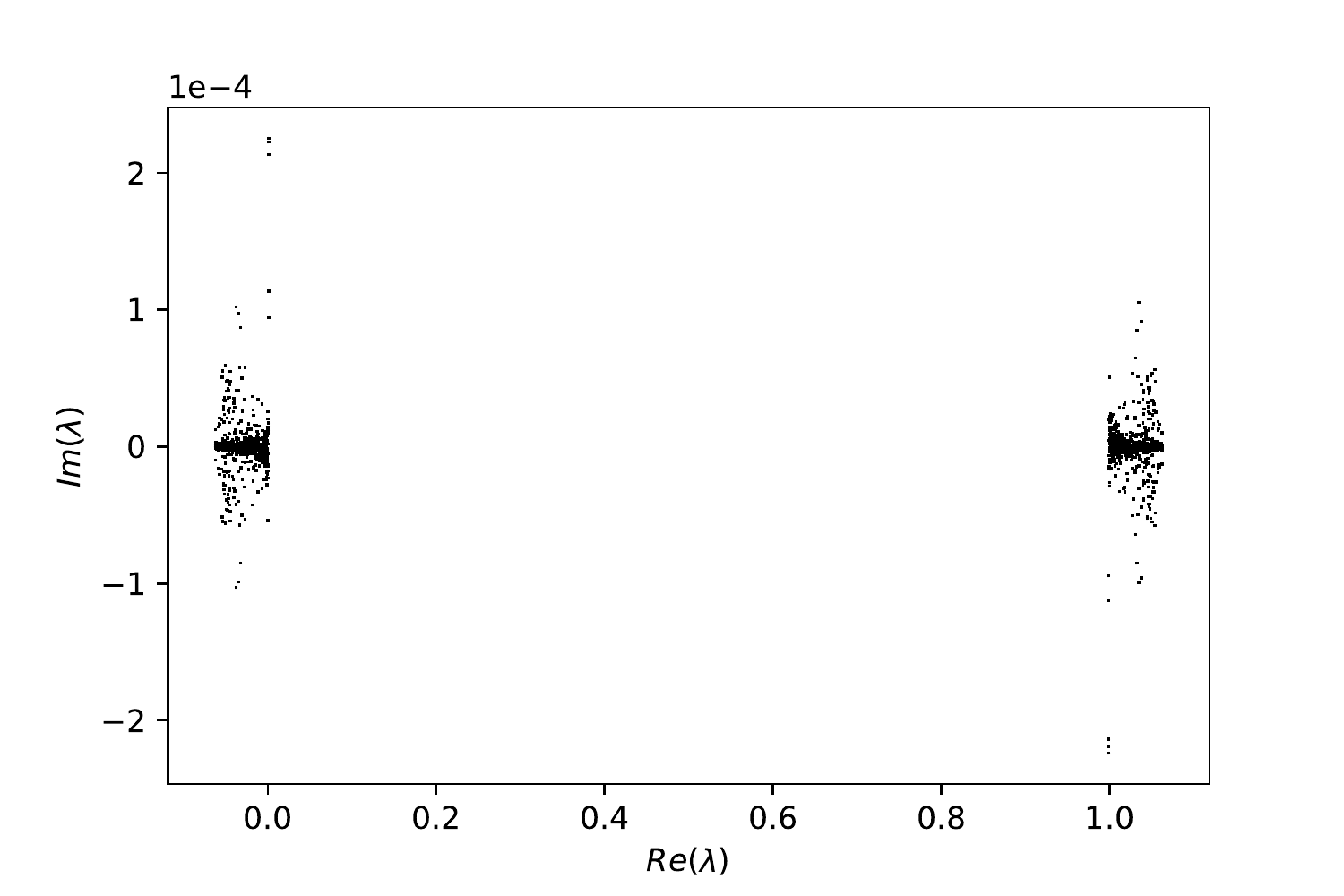}
\caption{Top: Singular values of the discrete strong form of the Calder\'on projector on the unit cube. Bottom: Eigenvalues of the discrete strong form.}
\label{fig:singeigvals}
\end{figure}

To demonstrate the action of the Calder\'on projector to a pair of non-compatible
Cauchy data we define two grid functions, both of which are constant one on the boundary.
\begin{lstlisting}[language=Python]
f1 = bempp.api.GridFunction.from_ones(
    calderon.domain_spaces[0])
f2 = bempp.api.GridFunction.from_ones(
    calderon.domain_spaces[1])
\end{lstlisting}
The two functions are defined on the pair of domain spaces discretising the product space
$H^{1/2}(\Gamma) \times H^{-1/2}(\Gamma)$. We can now apply the Calder\'{o}n projector to this pair of spaces to compute new grid functions which form a numerically compatible pair of Cauchy data for an exterior Helmholtz solution. The code snippet for this operation
is given by
\begin{lstlisting}[language=Python]
[u1, v1] = calderon * [f1, f2]
\end{lstlisting}
The grid functions \verb@u1@ and \verb@v1@ again live in the spaces of piecewise continuous and piecewise constant functions, respectively. We now apply the Calder\'on projector again to obtain
\begin{lstlisting}[language=Python]
[u2, v2] = calderon * [u1, v1]
\end{lstlisting}
The grid functions \verb@u1@ and \verb@u2@, respectively \verb@v1@ and \verb@v2@ should only differ in the order of the discretisation error. We can easily check this.
\begin{lstlisting}[language=Python]
error_dirichlet = (u2-u1).l2_norm() / u2.l2_norm()
error_neumann = (v2-v1).l2_norm() / v2.l2_norm()
\end{lstlisting}
For the corresponding values we obtain $1.2\times 10^{-4}$ and $8.0\times 10^{-4}$.
It is interesting to consider the singular values and eigenvalues of the discrete strong form of the Calder\'on projector. We can compute them easily as follows.
\begin{lstlisting}[language=Python]
from scipy.linalg import svdvals, eigvals
calderon_dense = bempp.api.as_matrix(calderon.strong_form())
sing_vals = svdvals(calderon_dense)
eig_vals = eigvals(calderon_dense)
\end{lstlisting}
The grid has 736 nodes. This means that the discrete basis for the possible Dirichlet data has dimension 736. For each
Dirichlet basis function there is a unique associated Neumann function via the Dirichlet-to-Neumann map. Hence, we expect the range of the Calder\'on projector to be of dimension 736 with all other singular values being close to the discretisation error. Correspondingly, for the eigenvalues we expect 736 eigenvalues close to 1 with all other eigenvalues being close to 0. This is indeed what happens as shown in Figure \ref{fig:singeigvals}. In the top plot we show the singular values of the discrete Calder\'on projector and in the right plot the eigenvalues. While the eigenvalues cluster around 1 and 0 the singular values show a significant drop-off between $\sigma_{736}\approx 1.04$ and $\sigma_{737}\approx 4.9\times 10^{-3}$, which corresponds to the approximation error as the accuracy of the hierarchical matrix approximation was chosen to be $10^{-3}$.

Finally, we would like to stress that while the eigenvalues of the discrete strong form are approximations to the eigenvalues of the continuous operator, the singular values of the discrete strong form are generally not. Given any operator $\bop{A}$ acting on a Hilbert space $\mathcal{H}$ the Galerkin approximation of the continuous eigenvalue problem $\bop{A}\phi = \lambda \phi$ is given as $\mat{A}x = \lambda M_Ax$, where $M_A$ is the mass matrix between the dual space and $\mathcal{H}$ with respect to the chosen dual form. If $M_A$ is invertible this is equivalent to $M_A^{-1}\mat{A}x = \lambda x$ or $\mat{A}^Sx=\lambda x$. The situation is more complicated for the singular values. For simplicity consider a compact operator (e.g. the single layer boundary operator) acting on $L^2(\Gamma)$. We have that
$$
\|A\|_{L^2(\Gamma)} = \sup_{\phi\in L^2(\Gamma)} \frac{\|A\phi\|_{L^2(\Gamma)}}{\|\phi\|_{L^2(\Gamma)}}.
 $$
 Let $M=C^TC$ be the Cholesky decomposition of the $L^2(\Gamma)$ mass matrix $M$ in a given discrete basis and $\mat{A}$ the Galerkin approximation in the same basis. Since $\|\phi\|_{L^2(\Gamma)} = \|C\bx\|_2$ for a function $\phi$ living in the discrete subspace of $L^2(\Gamma)$ with given coefficient vector $\bx$ it follows that
 $$
 \|A\|_{L^2(\Gamma)}\approx \max_{\bx\neq 0}\frac{\|CM^{-1}A\bx\|_2}{\|C\bx\|_2}=\|C^{-T}\mat{A}C^{-1}\|_2,
 $$
 which is generally not the same as $\|M^{-1}\mat{A}\|_2$. So while the strong form correctly represents spectral information it does not recover norm or similar singular value based approximations.

%We note that even though the largest eigenvalues are around one the largest singular values are a factor ten larger. This is a good numerical indicator that the exterior Calder\'{o}n projector is not normal (see \cite{Betcke14} for nonnormality of integral operators in acoustic scattering).

\subsection{Calder\'on preconditioning for acoustic transmission problems}
\label{sec:transmission}
As a final application we consider the Calder\'on preconditioned formulation of the following acoustic transmssion problem.
\begin{align}
-\Delta u\exterior - k^2 u\exterior = 0,~\text{ in } \Omega\exterior,\nonumber\\
-\Delta u\interior - n^2k^2 u\interior = 0,~\text{ in } \Omega,\nonumber\\
\gamma_0\interior u\interior = \gamma_0\exterior u\exterior +\gamma_0\exterior u^{inc},~\text{ on } \Gamma,\nonumber\\
\gamma_1\interior u\interior = \gamma_1\exterior u\exterior +\gamma_1\exterior u^{inc},~\text{ on } \Gamma,\nonumber\\
\lim_{|\bx|\to\infty}|\bx|\left(\frac{\partial}{\partial|\bx|}u\exterior (\bx) - iku\exterior(\bx)\right) = 0.
\end{align}
Here, $n = c^{+} / c^-$ is the ratio of the speed of sound $c^{+}$ in the surrounding medium to the speed of sound $c^{-}$ in the interior medium. The incident field is denoted by $u^{inc}$. The formulation that we present is based on \cite{Costabel85}. A generalized framework for scattering through composites is discussed in \cite{Claeys13}.
We denote by $V^{-} := \begin{bmatrix}\gamma_0^{-}u^{-} & \gamma_1^{-}u^{-}\end{bmatrix}^T$, $V^{+} := \begin{bmatrix}\gamma_0^{+}u^{+} & \gamma_1^{+}u^{+}\end{bmatrix}^T$, and $V^{inc} := \begin{bmatrix}\gamma_0^{+}u^{inc} & \gamma_1^{+}u^{inc}\end{bmatrix}^T$ the Cauchy data of $u^{-}$, $u^{+}$ and $u^{inc}$. Let $A^{+}$ be the multitrace operator associated with the wavenumber $k^{+}:=k$ and $A^{-}$ the multitrace operator associated with $k^{-}:=nk$ as defined in \eqref{eq:multitrace_operator}. From the Calder\'{o}n projector it now follows that
\begin{align}
\left(\frac{1}{2}I + A^{-}\right)V^{-} &= V^{-}\nonumber\\
\left(\frac{1}{2}I - A^{+}\right)V^{+} &= V^{+}
\end{align}
Together with the interface condition $V^{-} = V^{+} +V^{inc}$ we can derive from these relationships the formulation
\begin{equation}
\label{eq:single_trace}
\left(A^{-} + A^{+}\right)V^{+} = \left(\frac{1}{2}I - A^{-}\right)V^{inc}.
\end{equation}
This formulation is well defined for all wavenumbers $k > 0$ \cite{Costabel85}. Moreover, it admits a simple preconditioning strategy  \cite{Claeys13} based on properties of the Calder\'{o}n projector as follows. We note that $A^{+}$ is a compact perturbation of $A^{-}$ \cite{Sau11}. We hence obtain
$$
\left(A^{-} + A^{+}\right)^2 = \left(A^{-} + compact\right)^2 = \frac{1}{4}I + compact.
$$
We can therefore precondition \eqref{eq:single_trace} by squaring the left-hand side to arrive at
\begin{equation}
\label{eq:prec_single_trace}
\left(A^{-} + A^{+}\right)^2V^{+} = \left(A^{-} + A^{+}\right)\left(\frac{1}{2}I - A^{-}\right)V^{inc}.
\end{equation}

With the block operator algebra in place in Bempp the main code snippet becomes
\begin{lstlisting}[language=Python]
A_minus = multitrace_operator(grid, n * k, spaces='dual')
A_plus = multitrace_operator(grid, k, spaces='dual')
ident = multitrace_identity(grid, spaces='dual')
op = A_minus + A_plus
rhs_op = op * (.5 * ident - A_minus)
sol, info = bempp.linalg.gmres(op * op, rhs_op * v_inc,
                            use_strong_form=True)
\end{lstlisting}
As in the single-operator case we can intuitively write the underlying
equations and solve them. All mass matrix transformations are being
taken care off automatically.
\begin{figure}
  \center
  \includegraphics[width=8cm]{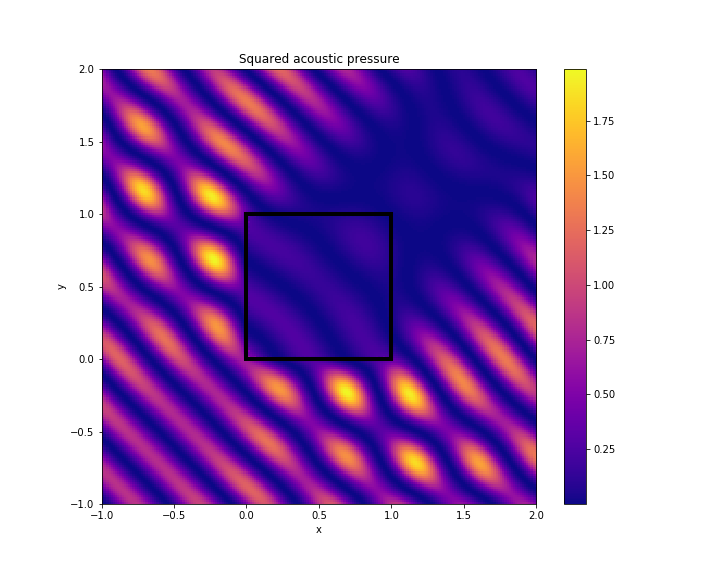}
  \caption{Squared acoustic pressure distribution of a wave
    travelling through a piecewise homogeneous medium.}
  \label{fig:square_scatt}
\end{figure}
An example is shown in Figure \ref{fig:square_scatt}. It demonstrates
a two-dimensional slice at height $0.5$ of a plane wave travelling
through the unit cube. In this example $k=10$ and $n=0.8$. The system
was solved in 7 GMRES iterations to a tolerance of $10^{-5}$.

\section{Conclusions}
\label{sec:conclusions}

In this paper we have demonstrated how a Galerkin based product
algebra can be defined and implemented. The underlying idea is very
simple. Instead of an operator being defined just trough a domain and a
test space we define it by a triplet of a domain space, range space,
and dual to range (test) space. This is more natural in terms of the
underlying mathematical description and allows the software
implementation of an automatic Galerkin operator product algebra.

We have demonstrated the power of this algebra using three examples,
the efficient evaluation of hypersingular boundary operators by
single-layer operators, the computation of the singular values and
eigenvalues of Calder\'{o}n projectors, and the Calder\'{o}n
preconditioned solution of an acoustic transmission problem.

As long as an efficient LU decomposition of the involved mass matrices
is possible  the product algebra can be implemented with little
overhead. Multiple LU decompositions of the same mass matrix can be
easily avoided through caching.

In this paper we focused on Galerkin discretizations of boundary
integral equations. Naturally, operator algebras are equally
applicable to Galerkin discretizations of partial differential
equations. The main difference here is that for large-scale three
dimensional problems an efficient LU decomposition of mass matrices
may not always be possible.

Finally, we would like to stress that the underlying principle of this
paper and its implementation in Bempp is to allow the user of software
libraries to work as closely to the mathematical formulation as
possible. Ideally, a user treats operators as continuous objects and
lets the software do the rest while the library ensures mathematical
correctness. The framework proposed in this paper and implemented in
Bempp provides a step towards this goal.

While in this paper we have focused on acoustic problems the extension
to Maxwell problems is straight forward and has been used in \cite{Scroggs2017}.

\bibliographystyle{siamplain}
\bibliography{references}

\begin{thebibliography}{10}

\bibitem{scipy}
{\em {SciPy}}.
\newblock \url{www.scipy.org}.

\bibitem{Buffa2007}
{\sc A.~Buffa and S.~H. Christiansen}, {\em A dual finite element complex on
  the barycentric refinement}, Mathematics of Computation, 76 (2007),
  pp.~1743--1769.

\bibitem{Claeys13}
{\sc X.~Claeys and R.~Hiptmair}, {\em Multi-trace boundary integral formulation
  for acoustic scattering by composite structures}, Communications on Pure and
  Applied Mathematics, 66 (2013), pp.~1163--1201.

\bibitem{Costabel85}
{\sc M.~Costabel and E.~Stephan}, {\em A direct boundary integral equation
  method for transmission problems}, Journal of Mathematical Analysis and
  Applications, 106 (1985), pp.~367 -- 413.

\bibitem{Greengard1987}
{\sc L.~Greengard and V.~Rokhlin}, {\em A fast algorithm for particle
  simulations}, Journal of Computational Physics, 73 (1987), pp.~325 -- 348.

\bibitem{Hackbusch2015}
{\sc W.~Hackbusch}, {\em Hierarchical matrices: algorithms and analysis},
  vol.~49 of Springer Series in Computational Mathematics, Springer,
  Heidelberg, 2015.

\bibitem{Hamdi1981}
{\sc M.~A. Hamdi}, {\em Une formulation variationnelle par {\'e}quations
  int{\'e}grales pour la r{\'e}solution de l’{\'e}quation de helmholtz avec
  des conditions aux limites mixtes}, CR Acad. Sci. Paris, s{\'e}rie II, 292
  (1981), pp.~17--20.

\bibitem{Hiptmair2006}
{\sc R.~Hiptmair}, {\em {Operator Preconditioning}}, Computers {\&} Mathematics
  with Applications, 52 (2006), pp.~699--706.

\bibitem{Kirby2010}
{\sc R.~C. Kirby}, {\em {From Functional Analysis to Iterative Methods}}, SIAM
  Review, 52 (2010), pp.~269--293.

\bibitem{Nedelec1982}
{\sc J.~C. Nedelec}, {\em Integral equations with non integrable kernels},
  Integral Equations and Operator Theory, 5 (1982), pp.~562--572.

\bibitem{Sau11}
{\sc S.~A. Sauter and C.~Schwab}, {\em Boundary element methods}, vol.~39 of
  Springer Series in Computational Mathematics, Springer-Verlag, Berlin, 2011.
\newblock Translated and expanded from the 2004 German original.

\bibitem{Scroggs2017}
{\sc M.~W. Scroggs, T.~Betcke, E.~Burman, W.~\'{S}migaj, and E.~van~’t Wout},
  {\em Software frameworks for integral equations in electromagnetic scattering
  based on {C}alder\'{o}n identities}, Computers \& Mathematics with
  Applications,  (2017).

\bibitem{Smigaj2015}
{\sc W.~{\'{S}}migaj, T.~Betcke, S.~Arridge, J.~Phillips, and M.~Schweiger},
  {\em {Solving Boundary Integral Problems with BEM++}}, ACM Transactions on
  Mathematical Software, 41 (2015), pp.~1--40.

\bibitem{Steinbach08}
{\sc O.~Steinbach}, {\em Numerical approximation methods for elliptic boundary
  value problems}, Springer, New York, 2008.
\newblock Finite and boundary elements, Translated from the 2003 German
  original.

\bibitem{Steinbach98}
{\sc O.~Steinbach and W.~L. Wendland}, {\em {The construction of some efficient
  preconditioners in the boundary element method}}, Advances in Computational
  Mathematics, 9 (1998), pp.~191--216.

\end{thebibliography}
\end{document}